\documentclass{amsart}
\usepackage{a4wide,amssymb,euscript,stmaryrd,epsfig}

\newtheorem{tm}{Theorem}[section]
\newtheorem{lemma}[tm]{Lemma}
\newtheorem{prop}[tm]{Proposition}
\newtheorem{cor}[tm]{Corollary}

 \theoremstyle{definition}
 \newtheorem{definition}{Definition}
\newtheorem{exmp}{Example}

\newtheorem{alg}{Algorithm}

\DeclareMathOperator*{\supp}{supp}
\DeclareMathOperator*{\osc}{osc}
\DeclareMathOperator*{\diam}{diam}

\begin{document}
\title{Sampling theorems on bounded domains}
\author{Massimo Fornasier and Laura Gori}

\date{}
\maketitle

\begin{abstract}
This paper concerns with iterative schemes for the perfect reconstruction of functions belonging to multiresolution spaces on bounded manifolds from nonuniform sampling. The schemes have optimal complexity in the sense that the computational cost to achieve a certain fixed accuracy is proportional to the computed quantity.
Since the iterations converge uniformly, one can produce corresponding iterative integration schemes that allow to recover the integral of functions belonging to multiresolution spaces from nonuniform sampling.
We present also an error analysis and, in particular, we estimate the $L^2$-error which one produces in recovering smooth functions in $H^s$, but not necessarily in any multiresolution space, and their integrals from nonuniform sampling. Several uni- and bi-variate numerical examples are illustrated and discussed.  We also show that one can construct a rather large variety of multiresolution spaces on manifolds from certain refinable bases on the real line formed by so-called GP-functions. This class of functions that contains in particular B-splines has remarkable properties in terms of producing well-conditioned bases. The resulting multiresolution analyses are well-suited for the application of the iterative recovering of functions from nonuniform sampling.
\end{abstract}


\section{Introduction}

The reconstruction of a {\it signal} from {\it sparse} and {\it nonuniform} sampling data is a well known problem \cite{BF2} and mainly addressed to  {\it guessing} or {\it learning} some relatively small missing part of a function by using the information of the relevant known part.
The problem of learning from examples is a corner stone in information theory and artificial intelligence. Very recently a probability theory of what can be learned from a relatively large distribution of data according to a (unknown) probability measure has been formalized in the beautiful contribution by Cucker and Smale \cite{CSm}. Inspired by this work, subsequent papers, for example \cite{BCDDT}, illustrate methods to construct estimators of the probabilistically best solutions of the learning problem.

Besides the statistical learning, mathematical methods and numerical algorithms based on deterministic techniques have been developed to compute missing parts of signals from few and sparse known sampling information and we refer, for example, to classical works of Feichtinger {\it et al.} \cite{FG92,FGS95,FS94,AF98} and Gr\"ochenig {\it et al.} \cite{AG01,GS}.
These methods and algorithms are essentially based on a {\it quasi-interpolation} of the function defined on the Euclidean space by means of suitable series expansions of irregularly shifted (translated) basic functions or, in its discrete version \cite{FGS95,BG}, by discrete finite series expansion of complex exponentials.
One of the typical applications is of course the restoration of digital signals and images, where some samples (ndr. pixels) can be corrupted, e.g., by scratches as it happens in old movies, or even completely missing, as in case of time-series of astronomical observation data. 
As a nice example of the use of such techniques for real applications, in \cite{Fo304} one of the authors discusses nonuniform sampling methods in combination with suitable variational models in order to restore colors of missing parts of destroyed art frescoes from those of few known fragments (previously detected and placed in their original sites \cite{FT03,FT04}), and the gray levels of the missing parts due to some pictures of the frescoes taken prior to the damage.
\\

In this paper we want to present sampling theorems for functions defined on compact subsets $\Omega$ of $\mathbb{R}^d$. Because of the boundedness nature of the domain such theorems are nicely suited for applications to signals in concrete situations. As we will show, even the proofs result relatively elementary with respect to those for functions defined on the whole Euclidean space \cite{FG92,AG01}.
The reconstruction from nonuniform sampling of functions defined on intervals has been also considered in \cite{FGS95,GS}.
In \cite{FGS95} the functions are modeled as complex trigonometric polynomials and in \cite{GS} as restrictions of functions defined on $\mathbb{R}$ as linear combination of shifted compactly supported functions. None of these two approaches can be really extended to the reconstruction of functions defined on manifolds, e.g., on the sphere. The approach that we present in this paper can be straightforward generalized to compact manifolds with or without boundary and it is computationally very efficient as the methods in \cite{FGS95,GS}.
In the following, the functions are assumed to belong to level spaces of a multiresolution analysis directly constructed on the domain (or manifold) \cite{CF1,CF2,DS0,DS2}, and not as a restriction of functions defined on larger spaces.
Moreover, the fact that one can associate wavelet bases to these multiresolution analyses to characterize function spaces of smooth functions, e.g., Sobolev spaces, allows to derive certain error estimations that would not be available otherwise.

The paper is organized as follows: Section 2 describes the abstract setting and the general requirements of the multiresolution spaces we use to model functions, and the corresponding sampling theorems for the reconstruction of functions from nonuniform sampling. Section 3 is devoted to the numerical implementation of the iterative reconstruction formulas derived in Section 2. Several numerical examples are shown and discussed. In Section 4 we show how the iterative formula can be also used to implement iterative numerical integration schemes from nonuniform sampling. We derive error estimators for smooth functions belonging to certain Sobolev spaces in Section 5. The construction on compact manifolds of multiresolution analyses with the properties that allow the application of the sampling theorems is addressed in Section 6. We will show that one can construct a rather large variety of multiresolution spaces from certain refinable functions on the real line chosen in the class described in \cite{GP00}.

\section*{Acknowledgement} 
M. Fornasier acknowledges the financial support provided through the
Intra-European Individual Marie Curie Fellowship Programme, 
under contract MEIF-CT-2004-501018, and the hospitality of NuHAG, Faculty of Mathematics, University of Vienna, Austria.

\section{Elementary axiomatic sampling theorems in $L^p(\Omega)$}

Let us assume that 

\begin{itemize} 
\item[MRA1)] is given a multiresolution analysis of finite dimensional spaces $V=\{V_j\}_{j \in \mathbb{N}}$ such that for all $p \in (1,\infty)$
$$
\chi_\Omega \in V_0(\Omega) \subset V_1(\Omega) \subset \dots V_j(\Omega) \subset L^p(\Omega), \quad \overline{\bigcup_j V_j(\Omega)}^{L^p(\Omega)} = L^p(\Omega);
$$
\item[MRA2)] each space has a suitable basis of nonnegative continuously differentiable functions $\Phi_j:=\{\phi_{j,0}, \dots, \phi_{j,N_j}\} \subset C^1(\Omega,\mathbb{R}_+)$,  such that $\Phi_{j-1} = A_j \Phi_j$ where $A_j \in M_{N_{j-1} \times N_j}(\mathbb{R})$ is a suitable (stochastic) scaling matrix, $2^{-\frac{d j}{2}} \sum_{k=0}^{N_j} \phi_{j,k} \equiv 1$ for all $j \in \mathbb{N}$. Moreover, we assume that $\diam(\supp(\phi_{j,k})) \leq C 2^{-j}$, $\#\{k: \supp(\phi_{j,k}) \cap \supp(\phi_{j,k'}) \neq \emptyset \} \leq m$ for some $m \in  \mathbb{N}$, uniformly with respect to $j$, and $\| \nabla \phi_{j,k}\|_\infty \leq C 2^{j (d/2+1)}$, uniformly with respect to both $k,j$;
\item[MRA3)] $P_j : C(\Omega) \rightarrow V_j(\Omega)$ is a bounded linear projector, i.e., $P_j f = f$ for all $f \in V_j(\Omega)$, for all $j \in \mathbb{N}$. Moreover we assume that
\begin{equation}
\label{eq1}
                                    f = \lim_{j \rightarrow \infty} P_j f, \quad \text{for all } f \in L^p(\Omega),
\end{equation}
with convergence in $L^p(\Omega)$.
\end{itemize}

\begin{definition}
A set $X=\{x_\ell\}_{\ell=0}^M$ of sampling nodes in $\Omega$ is called $\Delta$-dense if 
\begin{equation}
\label{eq2}
	\Omega = \left ( \bigcup_{\ell=0}^M  B_\delta(x_\ell) \right )\bigcap \Omega,
\end{equation} 
for $0<\delta:= \Delta^{-1}$, where $B_\delta(x)$ is the ball centered at $x$ of radius $\delta$.
\end{definition}
One has the following immediate result.

\begin{lemma}
For a given $\Delta$-dense set $X$ of sampling nodes, there exists a set of functions $\Psi=\{\psi_\ell\}_{\ell=0}^M \subset L^\infty(\Omega)$ with the following properties 
\begin{itemize}
\item[(1)] $0\leq \psi_\ell \leq 1$ for all $\ell=0,...,M$;
\item[(2)] $\supp(\psi_\ell) \subset B_\delta(x_\ell)$;
\item[(3)] $\sum_{\ell=0}^M \psi_\ell \equiv 1$.
\end{itemize}
\end{lemma}

Associated to a $\Delta$-dense set $X$ of sampling nodes one can define the following quasi-interpolation operator given by
\begin{equation}
\label{eq3}
Q_{\Psi,X} f := \sum_{\ell =0}^M f(x_\ell) \psi_\ell, \quad \text{for all } f \in C(\Omega).
\end{equation}
In the following we will make use also of another quasi-interpolation operator. For any $j \in \mathbb{N}$ and for all $k =0, ..., N_j$, let us choose a sampling point $\xi_{j,k} \in \supp(\phi_{j,k})$. Therefore, associated to the sampling set $\Xi_j=\{\xi_{j,k}\}_{k=0}^{N_j}$ we consider the following operator
\begin{equation}
\label{eq4}
S_{\Phi,\Xi_j} f := 2^{-d j/2} \sum_{\ell =k}^{N_j} f(\xi_{j,k}) \phi_{j,k}, \quad \text{for all } f \in C(\Omega).
\end{equation}

Let us formulate now the main result of this paper.

\begin{tm}
\label{tm1}
If $V=\{V_j\}_{j \in \mathbb{N}}$ is a multiresolution analysis in $L^p(\Omega)$ for $p \in (1,\infty)$ with properties MRA1)-MRA3), then for any $j \in \mathbb{N}$ there exists $\Delta_j>0$ large enough such that for any $\Delta_j$-dense sampling set $X_j:=X(\Delta_j)$ any function $f_j \in V_j(\Omega)$ can be exactly recovered from its samples $\{f(x_{j,\ell})\}_{\ell=1}^{M_j}$ by the following iterative algorithm
\begin{equation}
\label{eq5}
f^{(n+1)}_j = P_j Q_{\Psi_j,X_j} (f_j - f^{(n)}_j) + f^{(n)}_j, \quad n\geq 1, \quad f^{(0)}_j = P_j Q_{\Psi_j,X_j} f_j.
\end{equation} 
In fact one has
\begin{equation}
\label{eq6}
f_j := \lim_{n \rightarrow \infty} f^{(n)}_j,
\end{equation} 
where the convergence is uniform on $\Omega$, and then it is valid in $L^q(\Omega)$ for all $q \in [1,\infty]$. Moreover, for all $f \in L^p(\Omega)$, let $f_j =P_j f$, one has
\begin{equation}
\label{eq7}
f = \lim_{j \rightarrow \infty} \lim_{n \rightarrow \infty} f^{(n)}_j, 
\end{equation}
being $f^{(n)}_j$ defined by the iterative scheme \eqref{eq5} by means of suitable $\Delta_j$-dense sampling values ${\{f_j(x_{j,\ell})\}_{\ell=0}^{M_j}}$.
\end{tm}
The iterative scheme \eqref{eq5}  is analogous to that appearing in \cite{AF98}, but here it has been adapted for the reconstruction of functions defined on bounded domains.
\\

\begin{rem}
For $f \in C(\Omega)$ but $f \notin V_j(\Omega)$ for any $j$, the iterations \eqref{eq5} applied with initial value $ f^{(0)}_j = P_j Q_{\Psi_j,X_j} f$ will converge anyway to a function $f_j^{(\infty)}$ with the property 
$$P_j ( Q_{\Psi_j,X_j} f_j^{(\infty)} -   Q_{\Psi_j,X_j} f)=0.$$
\end{rem}

We shall now formulate in the following a reconstruction algorithm whose convergence is in fact a corollary of Theorem \ref{tm1}, and it is based on sampling at an {\it almost regular} (i.e., it can be understood as a perturbation of a uniform/regular) set of nodes.

\begin{cor}
\label{tm2}
If $V=\{V_j\}_{j \in \mathbb{N}}$ is a multiresolution analysis in $L^p(\Omega)$ for $p \in (1,\infty)$ with properties MRA1)-MRA3), then for any $j \in \mathbb{N}$ there exists a $j'\geq j$ large enough such that for any sampling set of the type $\Xi_{j'}=\{\xi_{j',k}\}_{k=0}^{N_{j'}}$ any function  $f_j \in V_j(\Omega)$ can be exactly recovered from its samples $\{f(\xi_{j',k})\}_{k=0}^{N_{j'}}$ by the following iterative algorithm
\begin{equation}
\label{eq8}
f^{(n+1)}_j =  P_j S_{\Phi_{j'},\Xi_{j'}}(f_j - f^{(n)}_j) + f^{(n)}_j, \quad n\geq 1, \quad f^{(0)}_j = P_j S_{\Phi_{j'},\Xi_{j'}} f_j.
\end{equation} 
In fact one has
\begin{equation}
\label{eq9}
f_j := \lim_{n \rightarrow \infty} f^{(n)}_j,
\end{equation} 
where the convergence is uniform on $\Omega$, and then it is valid in $L^q(\Omega)$ for all $q \in [1,\infty]$. 
Moreover, for all $f \in L^p(\Omega)$, let $f_j =P_j f$, one has
\begin{equation}
\label{eq10}
f = \lim_{j \rightarrow \infty} \lim_{n \rightarrow \infty} f^{(n)}_j, 
\end{equation}
being $f^{(n)}_j$ defined by the iterative scheme \eqref{eq8} by means of suitable  sampling values ${\{f_j(\xi_{j',k})\}_{k=0}^{N_{j'}}}$.
\end{cor}


Before working out the proof of Thereon \ref{tm1} we need to show the following technical lemma.

\begin{lemma}
\label{l1}
Under the notations and the assumptions of Theorem \ref{tm1}, for any $f \in V_j(\Omega)$, $j \in \mathbb{N}$, and any $\delta>0$ we define the oscillation function by
\begin{equation}
\label{eq11}
\osc_\delta(f)(x) = \sup_{y \in B_\delta(x)\cap \Omega}| f(x) - f(y)|.
\end{equation}
Then, one has 
\begin{equation}
\label{eq12}
\| \osc_\delta(f) \|_\infty \leq C_j \delta 2^{j(d/2+1)} |\Omega| \|f\|_\infty.
\end{equation}
\end{lemma}
\begin{proof}
Let us assume $f \in V_j(\Omega)$. Therefore, $f(x) = \sum_{k=0}^{N_j} a_k(f) \phi_{j,k}(x)$ and, if $y \in B_\delta(x)$, then 
\begin{eqnarray*}
&& |f(x)-f(y)| \leq \sum_{k=0}^{N_j} |a_k(f)| |  \phi_{j,k}(x)-   \phi_{j,k}(y)| \\
&\leq& \|\nabla \phi_{j,k}\|_\infty |x-y| \sum_{k=0}^{N_j} |a_k(f)| \leq  C \delta 2^{j(d/2+1)} \sum_{k=0}^{N_j} |a_k(f)| \\
&\leq&  (1+N_j)^{1/2} \delta 2^{j(d/2+1)} \left (\sum_{k=0}^{N_j} |a_k(f)|^2 \right )^{1/2} \leq  C B(1+N_j)^{1/2} \delta 2^{j(d/2+1)} \| f\|_2 \\
&\leq & C_j \delta 2^{j(d/2+1)} |\Omega| \|f\|_\infty.
\end{eqnarray*}
The estimation $\left (\sum_{k=0}^{N_j} |a_k(f)|^2 \right )^{1/2} \leq B \|f\|_2$ for some $B>0$ is valid because $\Phi_j$ is a basis for $V_j(\Omega)$. Here we have denoted $C_j = C B (1+N_j)^{1/2}$. Therefore, \eqref{eq12} is valid.
\end{proof}
Now we have all the necessary ingredients to prove Theorem \ref{tm1}.

\begin{proof}
We claim that for $\Delta_j>0$ large enough and for any $\Delta_j$-dense sampling set $X_j:=X(\Delta_j)$ one has
\begin{equation}
\label{eq13}
\| (I -  P_j Q_{\Psi_j,X_j}) f_j \| _\infty \leq \eta \|f_j \|_\infty, \quad \text{for all } f_j \in V_j(\Omega),
\end{equation}
for some $0<\eta<1$.
This would imply that $ P_j Q_{\Psi_j,X_j}$ is an invertible operator on $V_j(\Omega)$ with inverse
\begin{equation}
\label{eq14}
( P_j Q_{\Psi_j,X_j})^{-1} = \sum_{n=0}^\infty (I-  P_j Q_{\Psi_j,X_j})^n,
\end{equation}
so that
\begin{equation}
\label{eq15}
f_j = ( P_j Q_{\Psi_j,X_j})^{-1}  P_j Q_{\Psi_j,X_j} f_j = \sum_{n=0}^\infty (I-  P_j Q_{\Psi_j,X_j})^n P_j Q_{\Psi_j,X_j} f_j,
\end{equation}
with uniform convergence. Of course, it is immediate to see that \eqref{eq5} and \eqref{eq6} are just a reformulation of \eqref{eq15}.
So, let us show \eqref{eq13} and \eqref{eq14}:
\begin{eqnarray*}
\| (I -  P_j Q_{\Psi_j,X_j}) f_j \|_\infty &=& \| P_j (I - Q_{\Psi_j,X_j}) f_j \|_\infty\\
&\leq & C_j'' \| (I - Q_{\Psi_j,X_j}) f_j \|_\infty\\
& = & C_j'' \| f_j \left ( \sum_{\ell=0}^{M_j} \psi_\ell \right ) -\sum_{\ell=0}^{M_j} f_j(x_{j,\ell}) \psi_\ell\|_\infty \\
&\leq & C_j''  \| \sum_{\ell=0}^{M_j} |f_j(\cdot) -f_j(x_{j,\ell})| \psi_\ell\|_\infty\\
&\leq & C_j' \| \sum_{\ell=0}^{M_j} |\osc_{\delta_j} f_j(x_{j,\ell})| \psi_\ell\|_\infty \\
&\leq& C_j' \| \osc_{\delta_j} f_j\|_\infty.
\end{eqnarray*}

By applying Lemma \ref{l1}, one finally has
$$
\| (I -  P_j Q_{\Psi_j,X_j}) f_j \|_\infty \leq C_j \delta_j 2^{j(d/2+1)} |\Omega| \|f_j\|_\infty.
$$
Therefore, if $\Delta_j = \delta_j^{-1}>0$ is large enough, then  one immediately has  \eqref{eq14}.
The rest of the Theorem is shown as a consequence of MRA3).
Moreover, if $f \in C(\Omega)$, but not necessarily $f \in V_j(\Omega)$ for any $j$, then $ P_j Q_{\Psi_j,X_j} f \in V_j(\Omega)$ and 
\begin{eqnarray*}
P_j Q_{\Psi_j,X_j} f &=& \left( P_j Q_{\Psi_j,X_j} (P_j Q_{\Psi_j,X_j})^{-1} \right ) P_j Q_{\Psi_j,X_j} f\\
&=&  P_j Q_{\Psi_j,X_j} \left ( \sum_{n=0}^\infty (I-  P_j Q_{\Psi_j,X_j})^n P_j Q_{\Psi_j,X_j} f \right).
\end{eqnarray*}
Let us denote $f^{(\infty)}:= \left ( \sum_{n=0}^\infty (I-  P_j Q_{\Psi_j,X_j})^n P_j Q_{\Psi_j,X_j} f \right) \in V_j(\Omega)$. Then one has 
\begin{equation}
\label{eq16}
 P_j (  Q_{\Psi_j,X_j} f^{(\infty)} -  Q_{\Psi_j,X_j} f) =0.
\end{equation}
This proves the previous Remark.
\end{proof}

Similarly Corollary \ref{tm2} is shown:

\begin{proof}
Since we have assumed  $2^{-\frac{d j'}{2}} \sum_{k=0}^{N_{j'}} \phi_{j',k} \equiv 1$  one can write $\psi_k := 2^{-\frac{d j'}{2}}\phi_{j',k}$ and $x_{j',k}: = \xi_{j',k} \in \supp(\phi_{j',k})$, and one immediately see that $Q_{\Psi_j',X_j'} = S_{\Phi_{j'},\Xi_{j'}}$. Therefore an application of Theorem \ref{tm1} would conclude the proof. Let us anyway observe more explicitly that for $j' \geq j$ large enough, under the hypothesis MRA2) it is
$$
| f_j(x) -f_j(\xi_{j',k})|\phi_{j',k}(x) \neq 0,
$$
only if $x \in \supp(\phi_{j',k})$. That means that 
$$
| f_j(x) -f_j(\xi_{j',k})|\phi_{j',k}(x) \leq | \osc_{C m 2^{-j'}} f_j(\xi_{j',k})| \phi_{j',k}(x),
$$
pointwise. Therefore, following the proof of Theorem \ref{tm1} one obtains 
$$
\| (I -  P_j Q_{\Psi_j,X_j}) f_j \|_\infty \leq C_j m 2^{-j'} 2^{j(d/2+1)} |\Omega| \|f_j\|_\infty,
$$
so that $ C_j m 2^{-j'} 2^{j(d/2+1)} |\Omega|<1$ for $j'\geq j$ large enough.
\end{proof}



\subsection{On the multivariate interpolation problem on domains}

An interesting interpretation of Theorem \ref{tm1} is the following: \\

If, for a sequence of data $\{y_\ell\}_{\ell=0}^M$, there exists $f \in V_j(\Omega)$ such that $f(x_\ell)=y_\ell$ and $\{x_\ell\}_{\ell=0}^{M}$ are dense enough, then $f$ is unique. In fact if there were two functions $f,g \in V_j(\Omega)$ such that $f(x_\ell)=y_\ell = g(x_\ell)$ for all $\ell$, but $f \neq g$, one would have $Q_{\Psi_j,X_j} f =  Q_{\Psi_j,X_j} g$ and the following absurd consequence
$$
 f = (P_j Q_{\Psi_j,X_j})^{-1} P_j Q_{\Psi_j,X_j} f = (P_j Q_{\Psi_j,X_j})^{-1} P_j Q_{\Psi_j,X_j} g = g.
$$
In other words, if the interpolation problem of the data $\{(x_\ell,y_\ell)\}_{\ell=0}^M$ is solvable in $V_j(\Omega)$ then it is uniquely solvable. Unfortunately there are not simple conditions on $\{(x_\ell,y_\ell)\}_{\ell=0}^M$ so that the interpolation problem is solvable in $V_j(\Omega)$. However, as a consequence of \eqref{eq16}, one has the following {\it quasi}-interpolation argument:

For any $\Delta$-dense set $X$ of sampling nodes in $\Omega$ it is always possible to construct a piecewise constant interpolation operator
\begin{equation}
\label{eq3}
V_{X} f := \sum_{\ell =0}^M f(x_\ell) \chi_{\Omega_\ell}, \quad \text{for all } f \in C(\Omega),
\end{equation}
where 
$$
\Omega_\ell := \{x \in \Omega: |x_\ell - x| < | x_i -x|, \quad \forall i \neq \ell\}, \quad \ell =0 ,..., M
$$
defines a Voronoi decomposition of $\Omega$. In particular, one has $V_{X_j} f (x_\ell) = f(x_\ell)$, for all $\ell=0,...,M$.
It is not difficult to see that substituting $Q_{\Psi_j,X_j}$ with $V_{X_j}$ in Theorem \ref{tm1}, the result will be valid again.
Moreover, if $f = \sum_{\ell=0}^M y_\ell \chi_{\Omega_\ell}$ then, formally, it is $V_{X} f(x) =f(x)$ pointwise. Therefore one has
$$
0 = P_j ( V_{X_j} f^{(\infty)}- V_{X_j} f) =  P_j ( V_{X_j} f^{(\infty)}- f),
$$
uniformly.

\section{Numerical implementation, examples, and results}

In this section we want to illustrate a rather efficient numerical implementation of the scheme in \eqref{eq5}. 

\begin{figure}[ht]
\hbox to \hsize {\hfill \epsfig{file=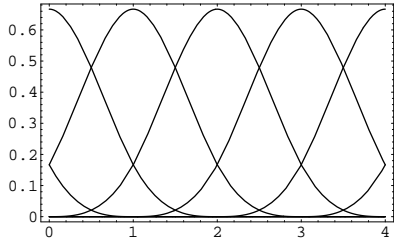,height=3.2cm} \hfill}
\caption{Cubic B-spline basis on the interval $[0,4]$}
\end{figure}

Without loss of generality we assume to work in the space $V_0(\Omega)$. First of all it is important to discuss how to construct a possible projector $P_0$. A natural choice can be in fact given by
\begin{equation}
\label{neq1}
P_0 f = \sum_{k=0}^{N_0} \langle f, \tilde \phi_{0,k} \rangle \phi_{0,k},  
\end{equation}
where $\tilde \Phi_0:=\{\tilde \phi_{0,k} \}_{k=0}^{N_0} \subset L^2(\Omega)$ is a biorthogonal dual basis for $\Phi_0$. In particular $P_0$ would result as a projection from $L^2(\Omega)$ onto $V_0(\Omega)$, of course bounded from $C(\Omega)$ into $V_0(\Omega)$ both being endowed with the sup-norm.

\begin{figure}[ht]
\hbox to \hsize {\hfill \epsfig{file=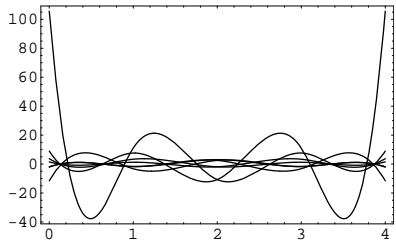,height=3.2cm} \hfill}
\caption{Canonical dual basis on the interval $\Omega= [0,4]$. The support of the functions coincides with the domain $\Omega$.}
\end{figure}

Therefore the problem of constructing a good projector is shifted now to the construction of suitable dual bases. For example the canonical dual biorthogonal basis in $V_0(\Omega)$ is computed as follows:

\begin{figure}[ht]
\hbox to \hsize {\hfill \epsfig{file=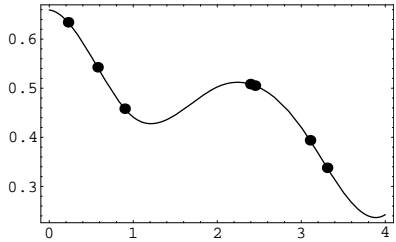,height=3.2cm} \hfill}
\caption{Sampling of a function in the cubic B-spline space $V_0[0,4]=\text{span}\{B_{0,3}(x+3-i):i=0,..,6\}$. Observe that the number of sampling nodes coincides with the dimension of the space $V_0$. However the nodes are strongly non-uniformly distributed with an almost coincidence of two of them.}
\end{figure}

Consider the Gramian matrix  $\mathbf{G}(\Phi_0)=(\langle \phi_{0,h}, \phi_{0,k} \rangle)_{h,k=0,...,N_0}$.
Then the rows (or the columns) of $\mathbf{G}(\Phi_0)^{-1}$ are the coordinates of the elements of the dual basis $\tilde \Phi_0$ with respect to $\Phi_0$, i.e.,
\begin{equation}
\label{neq2}
\tilde \Phi_0=\mathbf{G}(\Phi_0)^{-1} \Phi_0,
\end{equation}
where here, with an abuse of notation, we have considered both $\Phi_0$ and $\tilde \Phi_0$ as column vectors. 
Unfortunately the canonical dual basis can be formed by functions with support on the whole domain $\Omega$.
In order to obtain banded matrices and a more efficient reconstruction scheme, it will be more pleasant to deal with biorthogonal dual functions with support strictly contained in $\Omega$. The construction of such duals will be discussed in Section 6. 

\begin{figure}[ht]
\hbox to \hsize {\hfill \epsfig{file=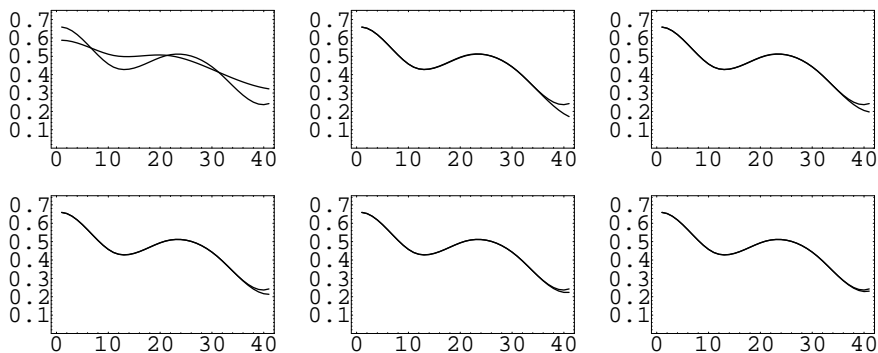,height=5cm} \hfill}
\caption{Successive iterations of the algorithm for the reconstruction of the function from its samples as illustrated in Figure 3. One can see the slow convergence of the approximant to the function to be restored.}
\end{figure}

Now that we have a recipe how to construct $P_0$ as in \eqref{neq1},  we should discuss how to implement numerically the general scheme \eqref{eq5}. Recalling from \eqref{eq3} that
$$
Q_{\Psi_0,X_0} f = \sum_{\ell =0}^{M_0} f(x_{0,\ell}) \psi_{0,\ell}, \quad \text{for all } f \in C(\Omega),
$$
one immediately has 
\begin{eqnarray*}
P_0 Q_{\Psi_0,X_0} f &=&  \sum_{k=0}^{N_0} \langle \sum_{\ell =0}^{M_0} f(x_{0,\ell}) \psi_{0,\ell}, \tilde \phi_{0,k} \rangle \phi_{0,k}\\
&=& \sum_{k=0}^{N_0} \left ( \sum_{\ell =0}^{M_0}f(x_{0,\ell}) \langle\psi_{0,\ell}, \tilde \phi_{0,k} \rangle \right) \phi_{0,k}.
\end{eqnarray*}
This suggests to define the notations $\mathbf{f}^s := [f(x_{0,0}) , ...,f(x_{0,M_0})]^T$, $\mathbf{M}_{\psi\tilde \phi}:=(\langle\psi_{0,\ell}, \tilde \phi_{0,k} \rangle)_{k=0,...,N_0;\ell=0,...,M_0}$, $\mathbf{\Phi_0}^s :=(\phi_{0,k}(x_{0,\ell}))_{\ell=0,...,M_0;k=0,...,N_0}$, and $\mathbf{\Phi_0}^c :=(\phi_{0,k}(\tau i))_{i=i \in \mathbb{Z}^d, \tau i \in \Omega;k=0,...,N_0}$ for $\tau>0$ small, and one can write
\begin{equation}
\label{neq3}
\mathbf{PQ}_s (\mathbf{f}^s) := \mathbf{\Phi_0}^s \mathbf{M}_{\psi\tilde \phi} \mathbf{f}^s ,
\end{equation}
\begin{equation}
\label{neq4}
\mathbf{PQ}_c (\mathbf{f}^s) :=  \mathbf{\Phi_0}^c \mathbf{M}_{\psi\tilde \phi} \mathbf{f}^s.
\end{equation}

If $f_0 \in V_0(\Omega)$ is the function that we want to reconstruct from its vector of samples $\mathbf{f_0^s} = [f_0(x_{0,0},...,f_0(x_{0,M_0})]^T$, one can use the following algorithm to compute/to approximate the vector $\mathbf{f_0^c} = [f_0(\tau i): i \in \mathbb{Z}^d, \tau i \in \Omega]^T$ of its sampling on a  $\tau^{-1}$-dense regular grid. The convergence of the following discrete scheme is of course ensured by the uniform and pointwise convergence of the original iterative algorithm \eqref{eq5}.

\begin{alg}
\label{alg1}
\noindent $\mathbf{RESTORE[n, \mathbf{f_0^s}] \rightarrow \mathbf{f_0^c}}$:\\
$\mathbf{f}^c:= \mathbf{PQ}_c (\mathbf{f_0^s});$\\
$\mathbf{f}^s:= \mathbf{PQ}_s (\mathbf{f_0^s});$\\
$i=0$;\\
While $i \leq n$ do\\
\indent $i:= i+1$\\
\indent $\mathbf{f}^c:= \mathbf{f}^c + \mathbf{PQ}_c (\mathbf{f_0^s} - \mathbf{f}^s)$;\\
\indent $\mathbf{f}^s:= \mathbf{f}^s + \mathbf{PQ}_s (\mathbf{f_0^s} - \mathbf{f}^s)$;\\
od \\
$\mathbf{f_0^c} := \mathbf{f}^c$.
\end{alg}

The discrete $\mathbf{PQ}$-procedure is implemented as in \eqref{neq3} and \eqref{neq4} by suitable matrix-matrix and matrix-vector multiplications. If the functions defining the matrices are assumed compactly supported (with support smaller than $\Omega$), then such matrices are banded.
Therefore, in such a case, Algorithm \ref{alg1} can be really implemented as a very fast reconstruction procedure. Unfortunately the canonical dual $\tilde \phi_{0,k} \in V_0(\Omega)$ is not in general locally supported and different biorthogonal duals to compute a projector $P_0$ should be considered.
Moreover, since the convergence is monotone with a rate of decay proportional to $\Delta^{-n}$, one has to expect of course faster results whenever $\Delta$-dense sets $X$ of sampling nodes are considered with $\Delta>0$ larger and larger. Figures 3-11 show uni- and bi-variate examples of applications of Algorithm \ref{alg1} in  cases where the sampling sets are highly nonuniform and not much dense, and, as a consequence, with a relatively slow convergence.

\begin{figure}[ht]
\hbox to \hsize {\hfill \epsfig{file=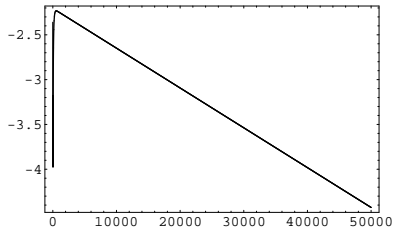,height=3.2cm} \hfill}
\caption{Maximal reconstruction error in logarithmic scale for the first 50000 iterations.}
\end{figure}

\begin{figure}[ht]
\hbox to \hsize {\hfill \epsfig{file=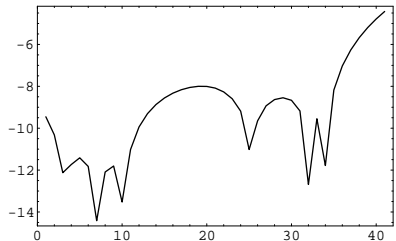,height=3.2cm} \hfill}
\caption{Pointwise reconstruction error in logarithmic scale at the $50000^{th}$ iteration.}
\end{figure}
\subsection{Computational cost}

Each iteration of the algorithm depends on certain matrix-vector multiplications where the involved matrices 
$$
\mathbf{PQ}_s (\mathbf{f}^s) := \mathbf{\Phi_0}^s \mathbf{M}_{\psi\tilde \phi} \mathbf{f}^s ,
$$
$$
\mathbf{PQ}_c (\mathbf{f}^s) :=  \mathbf{\Phi_0}^c \mathbf{M}_{\psi\tilde \phi} \mathbf{f}^s,
$$
can be both assumed to be banded (typically with bandwidth independent of the level $j$!) as soon as one uses suitable biorthogonal dual bases to define the projectors $P_j$. Moreover, the rate of convergence depends on the constant $\eta=\eta(\Delta)<1$ as in formula \eqref{eq13} so that
$$
\|f -  f^{(n)} \|_\infty \leq \eta^n \|f\|_\infty.
$$
Therefore, for a fixed density $\Delta$, for achieving a prescribed fixed tolerance $\varepsilon>0$ the computational cost can be assumed $O(\# supp(\mathbf{f}^c))$, where $\mathbf{f}^c$ is the output. This means that the algorithm has optimal complexity.

\begin{figure}[ht]
\hbox to \hsize {\hfill \epsfig{file=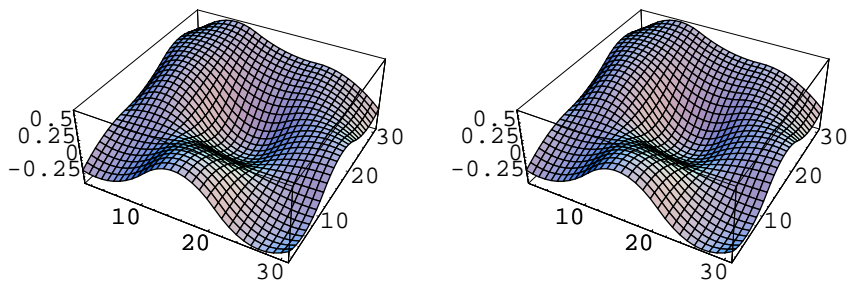,height=4.6cm} \hfill}
\caption{The original surface  belonging to $V_0([0,3]^2)=\text{span}\{B_{0,2}(x+2-i)B_{0,2}(y+2-j):i,j=0,..,4\}$ to be reconstructed  from its samples is shown on the left. On the right the result of the reconstruction is shown after 6000 iterations.}
\end{figure}

\begin{figure}[ht]
\hbox to \hsize {\hfill \epsfig{file=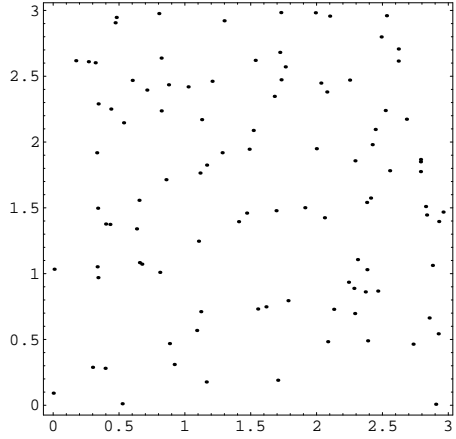,height=6cm} \hfill}
\caption{100=4 dim$(V_0)$ sparse random sampling points are chosen.}
\end{figure}

\begin{figure}[ht]
\hbox to \hsize {\hfill \epsfig{file=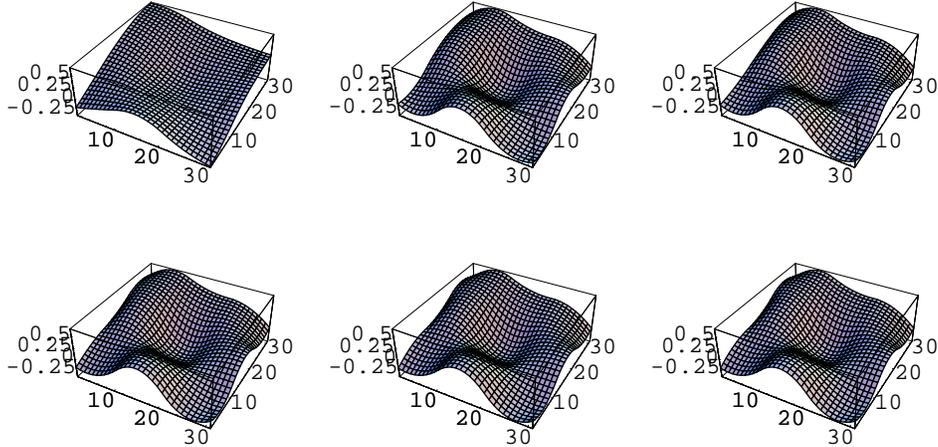,height=7cm} \hfill}
\caption{Some iterations of the reconstruction scheme.}
\end{figure}

\begin{figure}[ht]
\hbox  to \hsize {\hfill \epsfig{file=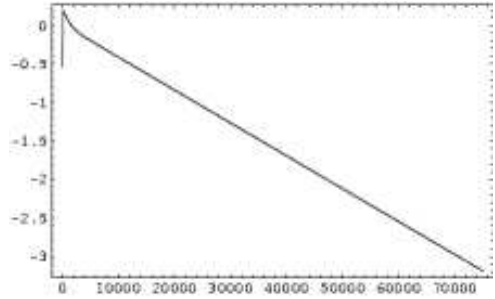, height=4cm}  \hfill }
\caption{The $\|\cdot\|_\infty$-error is shown in log- scale for successive iterations.}
\end{figure}

\begin{figure}[ht]
\hbox  to \hsize {\hfill \epsfig{file=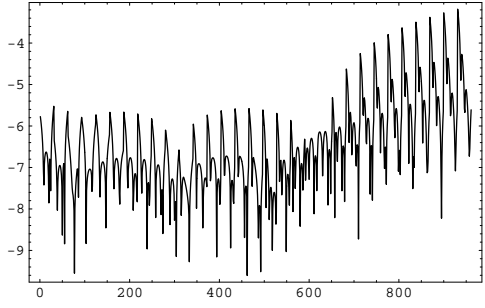, height=4cm}  \hfill }
\caption{The final local pointwise error in log-scale is shown.}
\end{figure}

\section{Multivariate iterative numerical integration}

Under the assumptions and notations considered so far, let $f_j \in V_j(\Omega)$. The iterative formula \eqref{eq5} in Theorem \ref{tm1} converges uniformly to $f_j$ and this implies also that 
\begin{equation}
\label{exact}
\int_\Omega f_j(x) dx = \lim_{n \rightarrow \infty} \int_\Omega f_j^{(n+1)}(x) dx. 
\end{equation}
On the other hand
\begin{eqnarray*}
 \int_\Omega f_j^{(n+1)}(x) dx &=& \int_{\Omega} P_j Q_{\Psi_j, X_j} (f_j - f_j^{(n)})(x) dx + \int_{\Omega} f_j^{(n)}(x) dx \\
&=& \sum_{k=0}^{N_j} \sum_{\ell=0}^{M_j} (f(x_{j,\ell}) - f_j^{(n)}(x_{j,\ell}) ) \langle \psi_{j,\ell}, \tilde \phi_{j,k} \rangle \left ( \int_\Omega \phi_{j,k}(x) dx \right ) + \int_{\Omega} f_j^{(n)}(x) dx\\
&=& \mathbf{\omega}_j^T \mathbf M_{\psi \tilde \phi}^j (\mathbf{f}_j^s - (\mathbf{f}_j^{(n)})^s) + \int_{\Omega} f_j^{(n)}(x) dx,
\end{eqnarray*}
where $\mathbf{M}_{\psi\tilde \phi}:=(\langle\psi_{j,\ell}, \tilde \phi_{j,k} \rangle)_{k=0,...,N_j;\ell=0,...,M_j}$, $\mathbf{f}_j^s := [f_j(x_{j,0}) , ...,f_j(x_{j,M_j})]^T$, $(\mathbf{f}_j^{(n)})^s := [f_j^{(n)}(x_{j,0}) , ...,f_j^{(n)}(x_{j,M_j})]^T$, and $\mathbf{\omega}_j = [ \int_\Omega \phi_{j,0}(x) dx, ...,  \int_\Omega \phi_{j,N_j}(x) dx]^T$. Therefore, similarly to Algorithm \ref{alg1}, one can formulate the following iterative integration procedure.

\begin{alg}
\label{alg2}
\noindent $\mathbf{INTEGRATE[n, \mathbf{f_j^s}] \rightarrow \mathbf{I}_j}$:\\
$\text{int}:= \mathbf{\omega}_j^T \mathbf{M}_{\psi \tilde \phi}^j \mathbf{f_j^s};$\\
$\mathbf{f}^s:= \mathbf{PQ}_s^j (\mathbf{f_j^s});$\\
$i=0$;\\
While $i \leq n$ do\\
\indent $i:= i+1$\\
\indent $\text{int} := \text{int} + \mathbf{\omega}_j^T \mathbf{M}_{\psi \tilde \phi}^j(\mathbf{f_j^s} - \mathbf{f}^s)$;\\
\indent $\mathbf{f}^s:= \mathbf{f}^s + \mathbf{PQ}_s^j (\mathbf{f_j^s}- \mathbf{f}^s)$;\\
od \\
$\mathbf{I}_j := \text{int}$.
\end{alg}

Here $\mathbf{PQ}_s^j$ is defined analogously to \eqref{neq3}.
\begin{prop}
\label{propINT}
Under the assumptions of Theorem \ref{tm1}, the procedure $\mathbf{INTEGRATE}$ has the following properties:
\begin{itemize}
\item[i)] It is linear, i.e., 
\begin{equation}
\label{linear}
\mathbf{INTEGRATE[n, \lambda \mathbf{f}^s + \mu \mathbf{g}^s] = \lambda \cdot INTEGRATE[n, \mathbf{f}^s] + \mu \cdot INTEGRATE[n,\mathbf{g}^s]},
\end{equation}
for all $\lambda, \mu \in \mathbb{C}$.
\item[ii)] for all $f_j \in V_j(\Omega)$ it is
\begin{equation}
\int_\Omega f_j(x) dx = \mathbf{INTEGRATE[\infty , \mathbf{f}_j^s]}:= \lim_{n \rightarrow \infty} \mathbf{INTEGRATE[n, \mathbf{f}_j^s]}.
\end{equation}
\end{itemize}
\end{prop}
\begin{proof}
To show i), it is sufficient to observe that all the operations executed in the procedure $\mathbf{INTEGRATE}$ are linear.
The proof of ii) is a straightforward application of \eqref{exact}. 
\end{proof}

\begin{rem}
The previous Proposition ii) means that integration formula $\mathbf{INTEGRATE[\infty , \mathbf{f}_j^s]}$ is exact on $f_j \in V_j(\Omega)$. Of course, one never can compute $\mathbf{INTEGRATE[\infty , \mathbf{f}_j^s]}$ but only $\mathbf{INTEGRATE[n, \mathbf{f}_j^s]}$ for $n$ large enough to ensure the desired accuracy. The accuracy is achieved without increasing the density of the sampling points $\mathbf{f}_j^s$ but just going further in the iterative integration process, see Table 1.
\end{rem}

\begin{cor}
\label{polyn}
If $\mathcal{P}_K(\Omega) \subset V_0(\Omega)$, where $\mathcal{P}_K(\Omega)$ denotes the set of algebraic polynomials of degree $K \in \mathbb{N}$ at most, then $\mathbf{INTEGRATE[\infty , \mathbf{p}^s]}$ is an exact iterative integration formula for all polynomials $p \in \mathcal{P}_K(\Omega)$.
\end{cor}

\begin{table}
	\begin{center}
		\begin{tabular}{|c|c|c|}
		\hline
	Iter.	& Approx. integral     & Error. \\
		\hline
	0 	&1.07874 &  0.131095\\		
		50	&0.97298 &  0.0253353\\	
100 &0.953824 &0.00617898\\
300 & 0.947611 & 0.0000340001\\
2000 &0.947643 &1.35891$\times 10^{-6}$\\
	\hline
		\end{tabular}
\end{center}
\begin{center}
Table 1. Application of the iterative integration. Integral value $\approx 0.9476446652462314$.
\end{center}

\end{table}

\section{Error analysis}

As we have shown in Section 2, for a continuous function $f \notin  V_j$ for all $j\geq 0$, Algorithm 1 applied on $\Delta_j>0$ dense samples of $f$ computes a function $f_j^{(\infty)}$ with the property  $P_j ( Q_{\Psi_j,X_j} f_j^{(\infty)} -   Q_{\Psi_j,X_j} f)=0$.  
This is an indirect information on the error that we have in the approximation $f \approx f_j^{(\infty)}$.
In this section we address the direct  estimation of the error $\|f - f_j^{(\infty)}\|_2$, under the assumption that $f$ has some regularity. In particular we will assume that $f \in H^s(\Omega)$, for $s \geq 0$, where $H^s(\Omega)$ denotes the Sobolev space of order $s$.  
Moreover, we also estimate the error that one produces by approximating the integral of $f$ by $\mathbf{INTEGRATE[\infty, \mathbf{f^s}]}$.
The key tool is the characterization of $H^s(\Omega)$ by suitable biorthogonal wavelets. In particular we assume in the following that\\
 
\begin{itemize}
\item[W1) ] there exists a wavelet biorthogonal basis $\{ \Psi_{j,k}\}_{j\geq 0,k \in \mathcal{J}_j}$ associated to the MRA $\{V_j(\Omega)\}_{j \geq 0}$, so that any function $f \in L^2(\Omega)$ can be written 
\begin{equation}
 f = \sum_{k=0,...,N_J} \langle f, \tilde \phi_{J,k} \rangle \phi_{J,k} + \sum_{j \geq J, k \in \mathcal{J}_j} \langle f, \tilde \Psi_{j,k} \rangle \Psi_{j,k} = P_J(f) + f_J^\perp;
\end{equation} 
\item[W2) ] $f \in H^s(\Omega)$, $s \geq 0$, if and only if 
\begin{equation}
c_1 \| f\|_{H^s(\Omega)} \leq \| P_J f\|_2+ \left( \sum_{j \geq J, k \in \mathcal{J}_j} 2^{2 s j} |\langle f, \tilde \Psi_{j,k} \rangle|^2 \right )^{1/2} \leq c_2 \| f\|_{H^s(\Omega)},
\end{equation}
for $c_1,c_2 >0$ independent of $f$.
This characterization of $H^s(\Omega)$ ensures in particular  that
\begin{eqnarray*}
\left( \sum_{j \geq J, k \in \mathcal{J}_j} |\langle f, \tilde \Psi_{j,k} \rangle|^2 \right )^{1/2} &=& 2^{-s J} \left( \sum_{j \geq J, k \in \mathcal{J}_j} 2^{2 s J} |\langle f, \tilde \Psi_{j,k} \rangle|^2 \right )^{1/2} \\
&\leq&  2^{-s J} \left( \sum_{j \geq J, k \in \mathcal{J}_j} 2^{2 s j} |\langle f, \tilde \Psi_{j,k} \rangle|^2 \right )^{1/2} \\
&\leq& c_2  2^{-s J} \| f\|_{H^s(\Omega)};
\end{eqnarray*}
\item[W3) ]$\| \Psi_{j,k}\|_\infty \leq C_\Psi 2^{\frac{d j}{2}}$ and $\# \mathcal{J}_j \leq C_{\mathcal{J}} 2^{d j}$;
\item[W4) ] $\int_\Omega \Psi_{j,k}(x) dx = 0$ for all $j\geq 0$ and $k \in \mathcal{J}_j$.
\end{itemize}

\begin{prop}
\label{estim}
Let $f \in C(\Omega) \cap H^s(\Omega)$ for $s\geq 0$, but not necessarily $f \in V_j$ for any $j \geq 0$. Under the assumptions of Theorem \ref{tm1} one has that $f_j^{(\infty)} := S_j f=\left ( \sum_{n=0}^\infty (I-  P_j Q_{\Psi_j,X_j})^n P_j Q_{\Psi_j,X_j} f \right) \in V_j(\Omega)$ has the property
\begin{equation}
\label{errest}
	\| f - f_j^{(\infty)} \|_2 \leq C \|S_j\| 2^{-s j}\|f\|_{H^s(\Omega)},
\end{equation}
where $\|S_j\|$ is the norm of the operator $S_j$ on $L^2(\Omega)$.
\end{prop}
\begin{proof}
With similar arguments as in the proofs of Lemma \ref{l1} and Theorem \ref{tm1}, one shows that the map $S_j:f \rightarrow f_j^{(\infty)}$ is bounded with the $L^2$-norm. Moreover this map coincides with the identity on $V_j$. This together with the properties W1-2 imply
\begin{eqnarray*}
\| f - f_j^{(\infty)}\|_2 &=& \| f - S_j (f -P_j f + P_j f)\|_2\\
&=& \|S_j(f -P_j f) - (f- S_j(P_j f))\|_2\\
&=&  \|S_j (f -P_j f) - (f- P_j f) \|_2\\
&\leq& (1+\|S_j\|) \| f- P_j f \|_2 \\
&\leq & C' (1+\|S_j\|) 2^{-s j} \| f\|_{H^s(\Omega)}.
\end{eqnarray*}
\end{proof}

\begin{lemma}
\label{estimINT}
Let $f \in C(\Omega)$. Under the assumptions of Theorem \ref{tm1} one has that
\begin{equation}
\left | \mathbf{INTEGRATE[\infty , \mathbf{f}^s]} \right | \leq \frac{|\Omega| \|P_j\|}{1-\eta} \|f\|_\infty,
\end{equation}
where $0<\eta<1$ is as in formula \eqref{eq13}.
\end{lemma}
\begin{proof}
By the remark in Section 2 and by \eqref{exact}, one has 
\begin{eqnarray*}
\left | \mathbf{INTEGRATE[\infty , \mathbf{f}^s]} \right | &=& \left | \int_\Omega \left ( \sum_{n=0}^\infty (I-  P_j Q_{\Psi_j,X_j})^n P_j Q_{\Psi_j,X_j} f(x) \right) dx \right |\\
&\leq & |\Omega| \left \| \sum_{n=0}^\infty (I-  P_j Q_{\Psi_j,X_j})^n P_j Q_{\Psi_j,X_j} f \right \|_\infty\\
&\leq& \frac{|\Omega|}{1 - \eta} \|  P_j Q_{\Psi_j,X_j} f \|_\infty\\
&\leq & \frac{|\Omega| \|P_j\| \|Q_{\Psi_j,X_j}\|}{1-\eta} \|f\|_\infty.
\end{eqnarray*}
One concludes observing that $ \|Q_{\Psi_j,X_j}\| =1$.
\end{proof}

Then one has the following error estimation result.

\begin{prop}
\label{errorAN}
Assume that 
 $\frac{3 d}{2}< s$. 
If the iterative integration formula is executed on a $\Delta_J$-dense set $X_J$, for $J \geq 0$ and  for $\Delta_J>0$ large enough,  then there exists a constant $C>0$ such that, for all $f \in H^s(\Omega)$, it is
\begin{eqnarray*}
\left | \mathbf{INTEGRATE[\infty , \mathbf{f}^s]} - \int_\Omega f(x) dx \right | &\leq & C \frac{|\Omega| \| P_J\|}{1-\eta}  \|f\|_{H^s} \frac{2^{(\frac{3 d}{2}-s)(J+1)}}{1- 2^{\frac{3d}{2}-s}}.
\end{eqnarray*}
\end{prop}
\begin{proof}
First of all observe that $s> \frac{d}{2}$ and then, by the Sobolev embedding theorem, it is $f \in C(\Omega)$ and it makes sense to consider its sampling.
By Proposition \ref{propINT} i) that ensures the linearity of $\mathbf{INTEGRATE}$ and by W1) 
{\small
\begin{eqnarray*}
&& \left | \mathbf{INTEGRATE[\infty , \mathbf{f}^s]} - \int_\Omega f(x) dx \right | \\
&=&  \left | \mathbf{INTEGRATE[\infty,P_J(f)^s]}+   \mathbf{INTEGRATE[\infty,\sum_{j \geq J,k} \langle f, \tilde \Psi_{j,k} \rangle \Psi_{j,k}^s]} -  \left (\int_\Omega P_J(f)(x) dx + \sum_{j \geq J,k} \langle f, \tilde \Psi_{j,k} \rangle  \int_{\Omega} \Psi_{j,k}(x) dx \right) \right |.
\end{eqnarray*}
}
In the last equality we could exchange the integral with the sum because the sum converges absolutely due to W2) and W3):
$$
\sum_{j \geq J,k} |\langle f, \tilde \Psi_{j,k} \rangle| |\Psi_{j,k}(x)| \leq C_\Psi \sum_{j \geq J,k} 2^{\frac{d}{j}}|\langle f, \tilde \Psi_{j,k} \rangle| \leq C \|f\|_{H^s} \sum_{j \geq J,k} 2^{(\frac{d}{2}-s)j} \leq  C' \|f\|_{H^s} \sum_{j \geq J} 2^{(\frac{3 d}{2}  -s)j}< \infty.
$$
Moreover, since $ \mathbf{INTEGRATE}$ is exact on $V_J$ by Proposition \ref{propINT} ii) and $\int_{\Omega} \Psi_{j,k}(x) dx=0$ one has
 \begin{eqnarray}
\label{last}
\left | \mathbf{INTEGRATE[\infty , \mathbf{f}^s]} - \int_\Omega f(x) dx \right | =  \left |   \mathbf{INTEGRATE[\infty,\sum_{j \geq J,k} \langle f, \tilde \Psi_{j,k} \rangle \Psi_{j,k}^s]}\right |
\end{eqnarray}
By Lemma \ref{estimINT}, equation \eqref{last}, and by using again W2) and W3), it is
\begin{eqnarray*}
 \left | \mathbf{INTEGRATE[\infty , \mathbf{f}^s]} - \int_\Omega f(x) dx \right | &\leq &  \frac{|\Omega|\| P_J\|}{1-\eta} \left \| \sum_{j \geq J,k} \langle f, \tilde \Psi_{j,k} \rangle \Psi_{j,k} \right \|_\infty\\
&\leq & C' \frac{|\Omega|\| P_J\|}{1-\eta}  \|f\|_{H^s} \sum_{j \geq J} 2^{(\frac{3 d}{2}-s)j}\\
&=&  C' \frac{|\Omega| \|P_J\|}{1-\eta}  \|f\|_{H^s} \frac{2^{(\frac{3 d}{2}-s)(J+1)}}{1- 2^{\frac{3d}{2}-s}}.
\end{eqnarray*}
\end{proof}

\begin{rem}
In other words, the previous theorem tells that if $f \in H^s(\Omega)$ with $s>\frac{3d}{2}$ and one has just a $\Delta_J$-dense set of sample points $\{f(x_{J,\ell})\}_{\ell=0}^{M_J}$, for $\Delta_J$ large enough, then its (multivariate) integral can be computed by using the procedure $\mathbf{INTEGRATE}$ with an error of order $O(\|P_J\| 2^{(\frac{3 d}{2}-s)(J+1)})$, once a sufficient number of iterations has been done.
\end{rem}

\section{Construction of MRA on domains and manifolds}

In this section we want to recall the construction of a rather large class of multiresolution analyses with properties MRA1-3 for which the axiomatic sampling theorems illustrated in the previous sections can be applied. It is also ensured that corresponding wavelet bases with properties W1-3 are also available.

\subsection{Composite multiresoultion analyses}

The main ingredient is the construction of biorthogonal refinable bases $\Phi_j^{Z_i}:=\{\phi_{j,0}, \dots, \phi_{j,N_j}\}$ and $\tilde \Phi_j^{\tilde Z_i}:=\{\tilde \phi_{j,0}, \dots, \phi_{j,N_j}\}$ to define the multiresolution spaces with complementary boundary conditions on $\Box=(0,1)^d$. In \cite{DS0} such bases have been derived from integer shifts of B-splines $\Theta_N$ on the real line. First one considers a biorthogonal dual bases $\tilde \Theta_{N, \tilde N}$ as described in \cite{CDF}, and then both $\Theta_N$ and $\tilde \Theta_{N, \tilde N}$ are restricted to define bases on $(0,1)$. Those elements of $\Theta_N$ and $\tilde \Theta_{N, \tilde N}$ that do intersect the boundary of $(0,1)$ are modified in order to ensure i) the polynomial reproduction with maximal degree, and ii) to match prescribed complementary boundary conditions. Such modifications of course imply that the resulting bases are no more biorthogonal. However such new bases can be in turn biorthogonalized. Finally, the bases $\Phi_j^{Z_i}$ and $\tilde \Phi_j^{\tilde Z_i}$ on $\Box$ are defined as tensor products of the resulting biorthogonal dual bases on $(0,1)$.
In \cite{CF1,CF2,DS0,DS2} the construction of spline-type multiresolution analyses on rather general manifolds has been proposed. Such construction is based on the decomposition of the manifold $\Omega$ into submanifolds $\Omega_i$ smoothly parametrized by $\Box:=(0,1)^d$, i.e., $\Omega = \cup_{i=1}^M \overline{\Omega_i}$, $\Omega_i = \kappa_i(\Box)$, $i=1, ..., M$, where $\Omega_i \cap \Omega_j = \emptyset$ for $i \neq j$, and $\kappa_i: \mathbb{R}^d \rightarrow \mathbb{R}^{d'}$, $d \leq d'$ are smooth regular functions.
The idea is to define biorthogonal  multiresolution analyses $V_j^{Z_i}(\Box)$ and $\tilde V_j^{\tilde Z_i}(\Box)$ on $\Box$ adjusted with suitable complementary boundary conditions $Z_i$ and $\tilde Z_i$ as described above \cite{DS0}, for each $i=1,...,M$. Each MRA is lifted on $\Omega_i$ by using the parametrization $\kappa_i$ to define $V_j(\Omega_i):=\{\varphi \circ \kappa_i^{-1}: \varphi \in V_j^{Z_i}(\Box)\}$ and similarly $\tilde V_j(\Omega_i)$. The boundary conditions are chosen to fit with a suitable bounded resolution of the identity $\mathcal{R} =\{R_i\}_{i=1}^M$ on $\Omega$ so that $V_j(\Omega):= \sum_{i=1}^M R_i V_j(\Omega_i)$ and $\tilde V_j(\Omega) := \sum_{i=1}^M R_i^* \tilde V_j(\Omega_i)$ will define global biorthogonal multiresolution analysis on $\Omega$ with properties MRA1-3. Moreover, associated to such global multiresolution analysis one can also define suitable biorthogonal wavelet bases \cite{DS2} with properties W1-3.
\\
Unfortunately the definition of the projectors $R_i$ makes use of certain extension operators that show increasing norms as soon as higher smoothness is required. This makes the resulting bases on $\Omega$ rather ill conditioned. One way to compensate this drawback is to start at the very beginning with refinable functions that are better conditioned than B-splines.

\begin{figure}[ht]
\hbox to \hsize {\hfill \epsfig{file=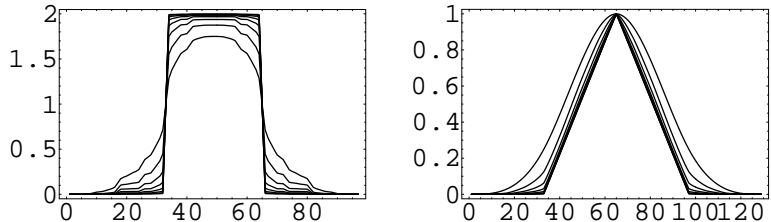,height=3.2cm} \hfill}
\caption{GP functions are illustrated for $n=2$ and $n=3$ and for different choices of $h>n-1$.}
\end{figure}

\subsection{GP bases}

A more general class of positive compactly supported refinable functions on the real line has been characterized in \cite{GP00} as the solutions of the refinement equations
\begin{equation}
\label{GP}
\varphi(x) = \sum_{k=0}^{n+1} a_k^{n,h} \varphi( 2 x- k),
\end{equation}
where
\begin{equation}
\label{GPmask}
a_{k}^{(n,h)} = 2^{-h} \left [   {n+1 \choose k} +4 (2^{h-n} -1) {n-1 \choose k-1}\right], \quad n \geq 2, h >n-1, k=0, ..., n+1.
\end{equation}
Let us call {\it GP-functions} the elements of this class. In particular B-splines are GP-functions for $h=n$.

\begin{figure}[ht]
\hbox to \hsize {\hfill \epsfig{file=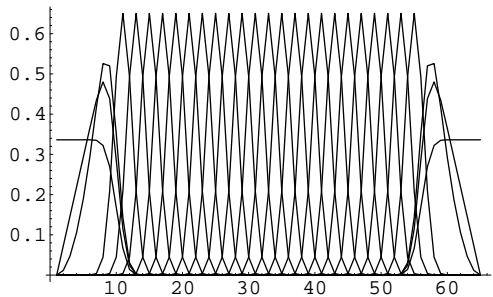,height=4cm} \hfill}
\caption{Normalized GP functions adapted to the interval  $(0,1)$ \cite[Section 4]{GPP04} for $n=5$ and $h=10$.}
\end{figure}
This generalization is worth for the following reason: It has been pointed out in \cite{GPP04} (see also \cite{GPP99,GPS01}) that they are better conditioned as bases $\Theta^{(n,h)}:=\{\varphi(x-k): k \in \mathbb{Z} \}$ for increasing values of $h$. This is due to the remarkable properties to have a smaller and smaller essential support, see Figure 12, while preserving the degree of smoothness, as soon as $h$ increases, see Table 2. 
\\

\begin{center}
\begin{tabular}{|c|c|c|c|c|c|c|c|} \hline\hline
{\em h} & 3&4&5&6&7&8&20 \\ \hline
$\kappa_2$ & 46.40 &26.06 &20.53 & 18.68 & 18.00 &17.73& 17.49 
\\ \hline\hline
\end{tabular}
\\

{Table 2. The spectral condition number $\kappa_2$ of the Gramian matrix for GP \\ functions adapted to the interval $(0,1)$ for $n=3$ and increasing values of $h$.}
\end{center}

In particular, in \cite{GPP04} it is shown how such bases can be adapted to the interval together with certain biorthogonal duals as derived in \cite{CL} although no complementary boundary conditions have been considered yet. We postpone this second adaptation to a successive contribution. \\
Let us shortly show instead a different method to produce biorthogonal duals for GP-bases as a generalization of the results in \cite{CDF}. As we have already stated in Section 3, the construction of a compactly supported biorthogonal dual is fundamental for the definition of suitable projectors $P_j$ raising banded matrices in Algorithms 1 and 2.

Associated to a mask $a = \{a_k: k=n_1,...,n_2\}$ for the refinement equation, one defines the symbol 
\begin{equation}
\label{symbol}
p(z):= \sum_{k=n_1}^{n_2} a_k z^k. 
\end{equation}
Equivalently one defines the trigonometric polynomial
\begin{equation}
\label{symbol2}
m(\xi) := \frac{1}{2} p(e^{- i \xi}).
\end{equation}
A refinable function $\tilde \varphi$ with mask $\tilde a^{(n,h)} = \{\tilde a_{k}^{(n,h)}: k=\tilde n_1,...,\tilde n_2\}$ defines a biorthogonal dual basis $\tilde \Theta^{(n,h)}:=\{\tilde \varphi(x-k): k \in \mathbb{Z} \}$, i.e., $\int_{\mathbb{R}} \varphi(x) \tilde \varphi(x-k) dx = \delta_{k,0}$, only if its associated symbol $\tilde m^{(n,h)}$ has the property
\begin{equation}
\label{maineq0}
 m^{(n,h)}(\xi) \overline{\tilde  m^{(n,h)}(\xi)} +  m^{(n,h)}(\xi + \pi) \overline{\tilde  m^{(n,h)}(\xi + \pi)} =1.
\end{equation}
Therefore, one can try to find a suitable trigonometric polynomial $\tilde m^{(n,h)}$ satisfying \eqref{maineq0}.
One has first the following result.
\begin{lemma}
\label{struct}
For $n>2$ and for $h>n-1$ it is 
\begin{equation}
\label{expl}
m^{(n,h)}(\xi) = e^{-i \frac{n+1}{2} \xi} \cos\left ( \frac{\xi}{2} \right)^{n-1} 1/(2^{h - n + 1})(2^{h - n + 1} - 1 + \cos(\xi)).
\end{equation}
\end{lemma}
\begin{proof}
In \cite[Formula (3.9)]{GP00} it is shown that 
$$
	p^{(n,h)}(z)=2^{-h}(1+z)^{n-1}(z^2+(2^{h-n+2}-2)z+1).
$$
A direct computation shows that $m^{(n,h)}(\xi) := \frac{1}{2} p^{(n,h)}(e^{- i \xi})$ is as in \eqref{expl}.
\end{proof}
Therefore, $m^{(n,h)}(\xi) = e^{-i \frac{n+1}{2} \xi} \cos\left ( \frac{\xi}{2} \right)^N r(\cos(x))$, where $N\in \mathbb{N}$ and $r$ is a certain polynomial. One can look for solutions of \eqref{maineq0} of the same form, i.e., $\tilde m^{(n,h)}(\xi) = e^{-i \frac{n+1}{2} \xi} \cos\left ( \frac{\xi}{2} \right)^{\tilde N} \tilde r(\cos(x))$, where $2 \ell:=N+\tilde N$ is even and $\tilde r$ is a suitable polynomial. Substituting these expressions in \eqref{maineq0} one obtains
\begin{equation}
\label{CDF}
\cos\left ( \frac{\xi}{2} \right )^{2 \ell} r(\cos(\xi)) \tilde r(\cos(\xi)) + \sin\left ( \frac{\xi}{2} \right )^{2 \ell} r(-\cos(\xi)) \tilde r(-\cos(\xi))=1.
\end{equation}
This equation is exactly the same as \cite[Formula (6.9)]{CDF}, and it is known to have solutions $\tilde r$ constrained by 
$$
r(\cos(\xi)) \tilde r(\cos(\xi)) = \sum_{i=0}^{\ell-1} {\ell -1 +i \choose i} \sin \left (\frac{\xi}{2} \right)^{2 i} + \sin \left (\frac{\xi}{2} \right)^{2 \ell} R(\cos(\xi)), 
$$
where $R$ is an odd polynomial. In other words one looks for a polynomial $\tilde r(x) = a_0 +a_1 x +...+ a_m x^m$ such that
$$
s[a_0,...,a_m](x) := r(x) ( a_0 +a_1 x +...+ a_m x^m) - \sum_{i=0}^{\ell-1} {\ell -1 +i \choose i} \left (\frac{1 - x}{2} \right)^{i}	
$$
is a polynomial divisible by $(\frac{1- x}{2})^\ell$ with an odd polynomial as the quotient. Since $r(x) = 1/(2^{h - n + 1})(2^{h - n + 1} - 1 + x)$ is a polynomial of degree 1 and the quotient $R$ has at least degree 1, it is clear that $\tilde r$ must be a polynomial of degree $m$ at least $\ell$. Let us assume then $m\geq \ell$. Therefore one reduces the problem to the solution of the following system of equations
\begin{equation}
\left \{ \begin{array}{ll}
\mod \left ( s[a_0,...,a_m](x), \left (\frac{1- x}{2} \right )^\ell \right) \equiv  0,\\ 
s[a_0,...,a_m](x)/\left (\frac{1- x}{2}\right)^\ell+ \quad s[a_0,...,a_m](-x)/\left (\frac{1- (-x)}{2}\right )^\ell  \equiv 0,\\
\end{array}
\right .
\end{equation}
where $p(x)/q(x)$ here indicates the quotient of the division of the polynomial $p$ by the polynomial $q$. Since the remainder $\mod ( s[a_0,...,a_m](x), (\frac{1- x}{2})^\ell)$ is a polynomial of degree $\ell-1$, the first equation will impose $\ell$ linear conditions on the $a_0,...,a_m$ variables. The quotient $s[a_0,...,a_m](x)/\left (\frac{1- x}{2}\right)^\ell$ is instead a polynomial of degree $m+1-\ell$. Therefore the second equation imposes further $\lfloor \frac{m+2-\ell}{2}\rfloor  +1$ equations. Altogether one generates $\lfloor \frac{m+1-\ell}{2}\rfloor  + \ell+1$ conditions. In particular for $m=\ell$, $m=\ell+1$, and $m=\ell+2$ one has that the number of conditions is exactly $m+1$, as many as the number of unknowns $a_0,...,a_m$. For $m \geq \ell+3$ the number of conditions is strictly less than $m+1$. This allows to impose further conditions to determine uniquely the solution.
Once determined a trigonometric polynomial $\tilde m^{(n,h)}$ by solving the linear equations, one should check that the equivalent conditions \cite[Theorem 4.3 C1-3]{CDF} are satisfied in order to conclude the successful construction of an admissible biorthogonal dual.

\begin{exmp}
Let us consider $n=3$ and $h=n+1=4$. By Lemma \ref{struct} it is
$$
	m^{(3,4)}(\xi) = e^{-2 i \xi} \cos(\xi/2)^2 \frac{3 + \cos(\xi)}{4}.
$$
Since in this case $N=n-1=2$, we choose $\tilde N =2$ and $2 \ell =N + \tilde N =4$.
Therefore one has to look for a polynomial $\tilde r(x) = a_0+a_1 x+ a_2 x^2$ of degree at least 2, such that
$$
s[a_0,a_1,a_2](x) = \frac{3 + x}{4} (a_0+a_1 x+ a_2 x^2) - \sum_{i=0}^{1} {i +1 \choose i} \left (\frac{1 - x}{2} \right)^{i}
$$
is divisible by $\left (\frac{1 - x}{2} \right)^{2}$, i.e., 
$$
\mod \left ( s[a_0,a_1,a_2](x), \left (\frac{1- x}{2}\right )^2 \right) = -2 +  3/4 a_0 - 1/4 a_1 - 5/4 a_2 + (1 + 1/4 a_0 + 5/4 a_1 + 9/4 a_2 ) x \equiv 0,
$$
and its remainder is an odd polynomial, i.e.,
$$
 s[a_0,a_1,a_2](x)/\left (\frac{1- x}{2}\right)^2+ \quad s[a_0,a_1,a_2](-x)/\left (\frac{1- (-x)}{2}\right )^2  = 2 a_2 + 10 a_3 \equiv 0
$$
These conditions are equivalent to the solution of the following linear problem:
$$
\left ( \begin{array}{ccc} 3/4 & -1/4 & -5/4 \\ 1/4 & 5/4 & 9/4 \\ 0 &2 &10 
\end{array} \right ) \left ( \begin{array}{c} a_0 \\ a_1 \\ a_2 \end{array} \right) =   \left ( \begin{array}{c} 2\\ -1 \\ 0 \end{array} \right).
$$
The solution is $(a_0,a_1,a_2)^T = (8/3, -25/12, 5/12)^T$. So that
$$
	\tilde m^{(3,4)} (\xi) = e^{-2 i \xi} \cos(\xi/2)^2 (8/3 -25/12 \cos(\xi)+5/12 \cos(\xi)^2).
$$
The corresponding mask is therefore given by
$$
	\tilde a^{(3,4)}=\{ a_{-1}^{(3,4)}=\frac{5}{96}, a_{0}^{(3,4)}=- \frac{5}{12} ,
   a_{1}^{(3,4)}=\frac{43}{96}, a_{2}^{(3,4)}=\frac{11}{6}, a_{3}^{(3,4)}=\frac{43}{96},
   a_{4}^{(3,4)}=-\frac{5}{12} , a_{5}^{(3,4)}=\frac{5}{96}\}
$$
that defines again a symmetric filter, see also \cite{CL}.
A direct numerical computation shows that $\tilde m^{(3,4)}$ so defined verify together with $m^{(3,4)}$ the conditions \cite[Theorem 4.3 C1-3]{CDF} so that 
$$
\tilde \varphi(x) = \sum_{k=-1}^{5} \tilde a_k^{3,4} \tilde \varphi( 2 x- k),
$$
is an $L^2(\mathbb{R})$ function and $\langle \varphi, \tilde \varphi(\cdot -k) \rangle = \delta_{k,0}$ for all $k \in \mathbb{Z}$.
\end{exmp}

As already mentioned above, the composite biorthogonal bases as derived on manifolds in \cite{DS2} might exhibit limited smoothness or high condition numbers. As proposed here, the use of better local bases certainly improves this drawback.
However an alternative way to overcome this difficulty is to use frames instead of biorthogonal bases, see for example \cite{DFR,S}. Such frames are constructed by Overlapping Domain Decompositions where the patches are again smooth images of $\Box$. The fact that one does not need to implement interfaces through patches to preserve global smoothness reduces the ill conditioning of the global system. The use of frames do not affect, e.g., the core of the proofs of Theorem \ref{tm1} and Proposition \ref{estim}.

\subsection{Implementing GP functions for the sampling problem}

As stated in the previous subsection GP-bases $\Theta^{(n,h)}$ adapted to the interval have better condition numbers for increasing values of $h>n-1$.
Thus one is tempted to affirm that for any $f \in V^{(n,h)}:=\text{span}\ \Theta^{(n,h)}$ for all $h>n-1$, $f$ will be reconstructed with increasing rate of convergence from a fixed $\Delta$-dense sampling set, for increasing values of $h$. In particular, it has been shown \cite{GP00} that $\mathcal{P}_{n-3} \subset   V^{(n,h)}$ for all $h>n-1$ and, therefore, one can test this claim on polynomials.
Surprisingly the claim is false. Numerical experiments for the reconstruction of the polynomial $f(x)=(x-\frac{1}{2})^2 -0.1$ on $(0,1)$ from a nonuniform sampling based on Algorithm 1 implemented by using $ V_j^{(5,h)}$, $j=5$, show that the asymptotic rate of convergence decreases for increasing values of $h$.
\begin{figure}[ht]
\hbox to \hsize {\hfill \epsfig{file=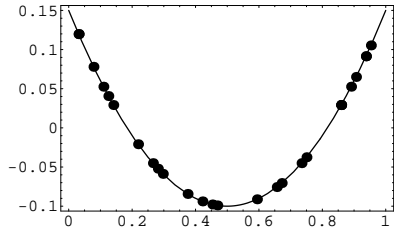,height=3.2cm} \hfill}
\caption{Nonuniform sampling of $f(x)=(x-\frac{1}{2})^2 -0.1$.}
\end{figure}
However, this also means that by using $V_j^{(5,h)}$, $j=5$, and $n-1<h<n$ one can anyway obtain better performances in the reconstruction than implementing Algorithm 1 by using just the B-spline space $V_j^{(5,5)}$, see Figure 15.

\begin{figure}[ht]
\hbox to \hsize {\hfill \epsfig{file=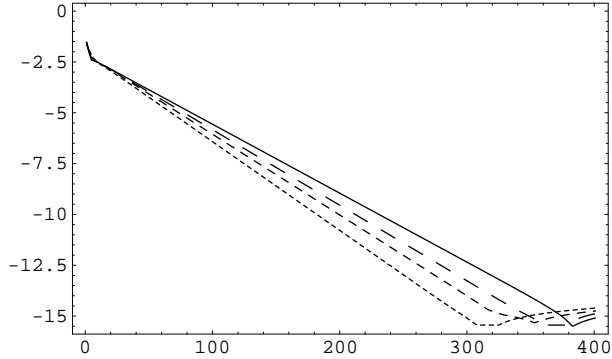,height=5cm} \hfill}
\caption{The solid line represents the uniform logarithmic error in the reconstruction of $f(x)=(x-\frac{1}{2})^2 -0.1$ by using the B-spline space $V_j^{(5,5)}$. The dashed lines represent  the uniform logarithmic error in the reconstruction of $f(x)=(x-\frac{1}{2})^2 -0.1$ by means of $V_j^{(5,h)}$ for $h \in \{4.1,4.3,4.5\}$. One can see that for $h=4.1$, one achieves the machine precision error with a 20$\%$ less of the number of iterations than for the choice $h=n=5$, corresponding to B-splines.}
\end{figure}

\subsection{Changing parameters in quasi-interpolation}

For a given sampling set $\{(x_\ell, y_\ell)\}_{\ell=0}^M$ and a fixed space $V_j$, it is not ensured that there exists $f_j \in V_j$ such that $f_j(x_\ell) = y_\ell$ for $\ell=0,...,M$.
Nevertheless, as soon as the set is dense enough Algorithm 1 will converge anyway and the resulting function $f_j^{(\infty)}$ will depend on several different parameters: The choice of the quasi-interpolation operator $Q_{\Psi,X}$ and the choice of the space $V_j$. 
Therefore one can tune the choice of $Q_{\Psi,X}$ and $V_j$ in order to optimize certain geometrical properties of the resulting $f_j^{(\infty)}$ with respect to the given data set $\{(x_\ell, y_\ell)\}_{\ell=0}^M$, see for example Figure 16.
\begin{figure}[ht]
\hbox to \hsize {\hfill \epsfig{file=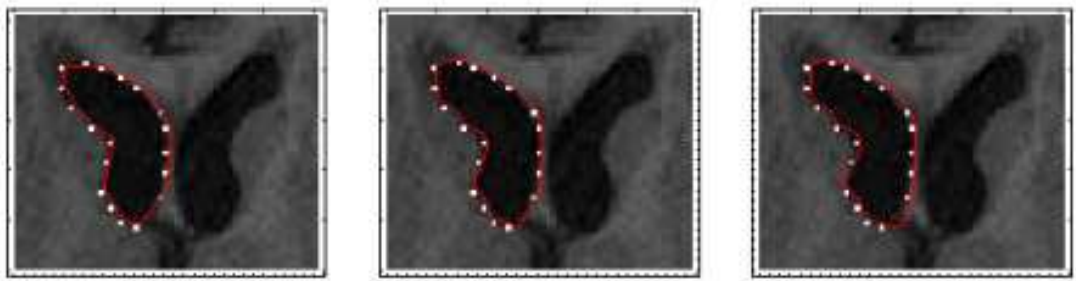,height=4cm} \hfill}
\caption{Segmentation of the brain phantom. The curves are computed with Algorithm 1 from the nonuniform sampling indicated by the white pixels. Different curves are generated by changing the GP space $V_j^{(n,h)}$ and the quasi-interpolation operator $Q_{\Psi,X}$.}
\end{figure}

$\vspace{1cm}$

\noindent Massimo Fornasier and Laura Gori\\
\noindent Dipartimento di Metodi e Modelli Matematici\\
per le Scienze Applicate\\
Universit\`a di Roma ``La Sapienza''\\
Via A. Scarpa 16/B - 00161 Roma\\
Italy
\\

\noindent Email: $\{${mfornasi@math.unipd.it, gori@dmmm.uniroma1.it}$\}$\\

\end{document}